\documentclass[12pt]{amsart}

\usepackage{amsthm, amsmath}
\usepackage{amssymb}
\usepackage[all]{xy}
\usepackage{mathrsfs}
\usepackage[latin1]{inputenc}
\usepackage{graphicx}
\usepackage{xspace}

\numberwithin{equation}{subsection}

\newcommand{\CC}{(\textbf{CC})\xspace}

\newcommand \bbq{{\mathsf {q}}}
\newcommand\CH{{\mathcal H}}

\theoremstyle{plain}

\newtheorem{theorem}{Theorem}[section]

\newtheorem{thm}[theorem]{Theorem}
\newtheorem{cor}[theorem]{Corollary}

\newtheorem{prop}[theorem]{Proposition}

\newtheorem{ques}[theorem]{Question}

\theoremstyle{definition}

\def\le{\leqslant}

\def\L{\Lambda}

\def\s{\sigma}

\def\th{\theta}
\def\k{\kappa}
\def\l{\lambda}

\def\ZZ{\mathbb Z}
\def\NN{\mathbb N}
\def\QQ{\mathbb Q}

\def\CC{\mathbb C}

\def\ci{\mathcal I}

\def\co{\mathcal O}

\def\tW{\tilde W}

\def\tS{\tilde S}

\usepackage{amssymb} % symboler fra AMS
\usepackage{latexsym} % LaTeX symboler
\usepackage{amsfonts} % fonte fra AMS
\usepackage{amsmath} % matematik fonte fra AMS
\usepackage{eucal} % 'krlle'-fonte (skriv \mathcal)
\usepackage{bm} % fed skrift i matematik
\usepackage{bbm} % skrift med dobbelt streg til f.eks. talomrï¿½er
\usepackage{graphicx} % mulighed for at inkludere grafik
\usepackage[english]{varioref} % smarte referencer - se nederst
\usepackage[nice]{nicefrac} % pï¿½ere brker
\usepackage[all]{xy}

\newcommand{\kk}{\Bbbk}

\def\<{\langle}
\def\>{\rangle}

\begin{document}

\title[The cocenter-representation duality]{The cocenter-representation duality}

\dedicatory{Dedicated to David Vogan on his 60th birthday}

\author[X. He]{Xuhua He}
\address{Department of Mathematics, University of Maryland, College Park, MD 20742, USA and Department of Mathematics, HKUST, Hong Kong}
\email{xuhuahe@gmail.com}

\begin{abstract}
The purpose of this note is to explain the relation between the cocenter and Grothendieck group of finite dimensional representations of affine Hecke algebras. It is based on my talk given in Vogan conference, 05/2014. 
\end{abstract}

%\date{\today}

\maketitle

\section*{Introduction}

Affine Hecke algebras arise naturally in the study of smooth representations of reductive $p$-adic groups. Finite dimensional complex representations of affine Hecke algebras (under some restriction on the isogeny class and the parameter function) has been studied by many mathematicians, including Kazhdan-Lusztig \cite{KL}, Ginzburg \cite{CG}, Lusztig \cite{L1}, Reeder \cite{Re}, Opdam-Solleveld \cite{OS}, Kato \cite{Kat}, etc. The approaches are either geometric or analytic. 

In this note, we'll discuss a different route, via the so-called ``cocenter-representation duality'', to study finite dimensional representations of affine Hecke algebras (for arbitrary isogeny class and for a generic  complex parameter). This route is more algebraic, and allows us to work with complex parameters, instead of equal parameters or positive parameters. We also expect that it can be eventually applied to the ``modular case'' (for representations over fields of positive characteristic and for parameter equal to a root of unity). 

\section{A naive example: group algebras of finite groups} 

\subsection{}\label{1.1} To explain some basic idea of the cocenter-representation duality, we start with a naive example. 

Let $G$ be a finite group and $V$ be a finite dimensional complex representation of $G$. We define the character of $V$ by $\chi_V(g)=Tr(g, V)$ for $g \in G$. Then $\chi_V$ is a class function, i.e., $\chi_V(g)=\chi_V(g')$ if $g$ and $g'$ are conjugate in $G$. 

Let $V_1, \cdots, V_k$ be the irreducible representations of $G$ (up to isomorphism) and $\co_1, \cdots, \co_l$ be the conjugacy classes of $G$. It is well known that $k=l$ and 

(a) The matrix $(\chi_{V_i}(g_j))_{1 \le i, j \le k}$ is invertible, where $g_j$ is a representative of $\co_j$.

This matrix is called the {\it character table} of $G$. 

\subsection{}\label{1.2} We reformulate \S\ref{1.1} (a) in a different way. 

Let $H_1=\CC[G]$ be the group algebra of $G$. Let $[H_1, H_1]$ be the commutator of $H_1$, the subspace of $H_1$ spanned by $[h, h']:=h h'-h' h$ for all $h, h' \in H_1$. We call the quotient space $\bar H_1=H_1/[H_1, H_1]$ the {\it cocenter} of $H_1$.  

It is easy to see that

(1) For any $g, g'$ in a given conjugacy class $\co$ of $G$, the image of $g$ and $g'$ in $\bar H_1$ are the same. We denote it by $[\co]$. 

(2) $\{[\co_1], \cdots, [\co_k]\}$ is a basis of $\bar H_1$. 

Let $R(H_1)=R(G)$ be the Grothendieck group of finite dimensional complex representations of $H_1$. Then $\{V_1, \cdots, V_k\}$ is a basis of $R(H_1)$. The trace map $$Tr: H_1 \to R(H_1)^*, \qquad g \mapsto (V \mapsto Tr(g, V))$$ factors through $Tr: \bar H_1 \to R(H_1)^*$. We reformulate \S\ref{1.1}(a) as 

(a) $Tr: \bar H_1 \to R(H_1)^*$ is an isomorphism of vector spaces. 

We call it the cocenter-representation (in short, $\bar H-R(H)$) duality of the group $G$ and the group algebra $H_1$. 

\section{Coxeter groups and Hecke algebras} 

\subsection{} Let $S$ be a finite set and $M=(m_{st})_{s, t \in S}$ be a Coxeter matrix, i.e., $m_{ss}=1$ for all $s \in S$ and $m_{st}=m_{ts} \in \{2, 3, \cdots\} \cup \{\infty\}$ for all $s \neq t$ in $S$. The Coxeter group $W=W(M)$ is the group defined by the generators $s \in S$ and relations $(s t)^{m_{st}}=1$ for $s, t \in S$ with $m_{st}<\infty$. 

In additional to the group structure, the length function $\ell: W \to \NN$ plays a crucial in this note. Here for any $w \in W$, $\ell(w)$ is the minimal integer $n$ such that $w=s_1 \cdots s_n$, where $s_i \in S$. 

\subsection{}\label{2.2} Fix a set of indeterminates $\bbq=\{\bbq(s): s\in S\}$
such that $\bbq(s)=\bbq(t)$ if $s$ and $t$ are conjugate in $W$, and let
$\Lambda=\ZZ[\bbq(s)^{\pm 1}: s\in S].$

The generic Hecke algebra
  $\CH=\CH(W,\bbq)$ is the $\Lambda$-algebra generated by $\{T_w:
  w\in \tW\}$ subject to the relations:
\begin{enumerate}
\item $T_w\cdot T_{w'}=T_{ww'}$, if $\ell(ww')=\ell(w)+\ell(w')$;
\item $(T_s+1)(T_s-\bbq(s)^2)=0,$ $s\in S$.
\end{enumerate}

%Here $\{T_w\}_{w \in W}$ is a basis of $\CH$. 

If we assign a nonzero element $c_s \in \CC^\times$ to $\bbq(s)$ for $s \in S$, then we can regard $\CC$ as a $\L$-module and $H_c:=\CH \otimes_\L \CC$ is the specialization of $\CH$. In particular, if $c_s=1$ for all $s \in S$, then we obtain the group algebra $H_1=\CC[W]$. 

\section{The $\bar H-R(H)$ duality for finite Hecke algebras} 

\subsection{} In this section, we assume that $W$ is a finite Coxeter group and $\CH$ is a finite Hecke algebra. 

By Tits' deformation theorem \cite[Theorem 7.4.6 \& 7.4.7]{GP00}, for generic parameter $\mathbf c=\{c_s; s \in S\}$, the Hecke algebra $H_{\mathbf c}$ is isomorphic to the group algebra $H_1$, and hence the number of irreducible representations of $H_{\mathbf c}$ equals the number of irreducible representations of $H_1$, and hence equals the number of conjugacy classes of $W$. 

However, the trace function $Tr(-, V)$ is not a ``class function'' in the sense that $Tr(T_w, V)=Tr(T_{w'}, V)$ if $w$ and $w'$ are conjugate in $W$. 

\

To overcome this difficulty, we use the following remarkable property of finite Weyl groups, first discovered by Geck and Pfeiffer \cite{GP93} via a case-by-case analysis, with the aid of computer for exceptional types. A case-free proof was found 20 years later in \cite{HN12}. 

\begin{thm}\label{f-min}
Let $W$ be a finite Coxeter group and $\co$ be a conjugacy class. Let $\co_{\min}$ be the set of minimal length elements in $\co$. Then 

(1) Each element $w$ of $\co$ can be brought to a minimal length element $w'$ by conjugation by simple reflections which reduce or keep constant the length. 

(2) Any two elements in $\co_{\min}$ are ``strongly conjugate'', i.e., conjugate in the associated Braid group. 
\end{thm}

\subsection{} In the rest of this section, we fix a generic parameter $\mathbf c$ and simply write $H$ for $H_{\mathbf c}$. 

By Theorem \ref{f-min} (2), if $w$ and $w'$ are of minimal length in a given conjugacy class $\co$ of $W$, then the image of $T_w$ and $T_{w'}$ in the cocenter $\bar H=H/[H, H]$ are the same. We denote it by $T_\co$. 

By Theorem \ref{f-min} (1), $\{T_\co\}$ spans $\bar H$. 

By comparing the number of irreducible representations of $H$ and the number of conjugacy classes of $W$, one deduces that $\{T_\co\}$ is in fact a basis of $\bar H$. Therefore, 

\begin{thm}
Let $W$ be a finite Coxeter group and $H_{\mathbf c}$ is the Hecke algebra of $W$ and generic parameter $\mathbf c$. Then $Tr: \bar H_{\mathbf c} \to R(H_{\mathbf c})^*$ is an isomorphism of vector spaces. 
\end{thm}

This leads to the definition and study of ``character table'' of finite Hecke algebras in \cite{GP93}. 

\section{Affine Hecke algebras} 

\subsection{} The main purpose of this note is to study the (extended) affine Hecke algebra of (extended) affine Weyl groups. Let us first recall the definition. 

Let $\Phi=(X, R, X^\vee, R^\vee, \Pi)$ be a based root datum and $W_0=W(R)$ be the finite Weyl group. The affine Weyl group and the extended affine Weyl group are defined to be $$W_a=\ZZ R \rtimes W_0, \qquad  \tW=X \rtimes W_0.$$

The affine Weyl group $W_a$ is a Coxeter group with a set of simple reflections $\tS$. The extended affine Weyl group $\tW$, strictly speaking, is not a Coxeter group. However, $\tW$ is a quasi-Coxeter group, in the sense that there is a natural length function on $\tW$ such that $\tW=W_a \rtimes \Omega$, where $\Omega$ is the set of length-zero elements in $\tW$. 

For example, the extended affine Weyl group associated to the $p$-adic group $GL_n(\mathbb Q_p)$ is the extended affine Weyl group $\ZZ^n \rtimes S_n$. This is the reason that we are interested in the extended affine Weyl groups, not just the affine Weyl groups. 

\subsection{} Fix a set of indeterminates $\bbq=\{\bbq(s): s\in \tS\}$
such that $\bbq(s)=\bbq(t)$ if $s$ and $t$ are conjugate in $\tW$, and let
$\Lambda=\ZZ[\bbq(s)^{\pm 1}: s\in \tS].$ We define the extended affine Hecke algebra $\tilde \CH$ and its specialization $\tilde H_{\mathbf c}$ in a similar way as in \S\ref{2.2}. The basis $\{T_w\}_{w \in \tW}$ gives the Iwahori-Matsumoto presentation of $\tilde \CH$. It is related to the quasi-Coxeter structure of $\tW$. 

For any $J \subset \tS$, we denote by $W_J$ the subgroup of $W_a$ generated by the simple reflections in $J$. If $W_J$ is finite, then we call it a parahoric subgroup of $\tW$ and we denote by $\CH_{J}^{fin}$ the parahoric subalgebra generated by $T_w$ for $w \in W_J$. 

\subsection{} Another presentation we need for $\CH$ is the Bernstein-Lusztig presentation.

(a) $\{\th_x T_w; x \in X, w \in W_0\}$ is a basis of $\CH$. 

The definition can be found in \cite{L1}. This presentation is related to the semi-product $\tW=X \rtimes W_0$. It plays an important role in the study of representations of affine Hecke algebras, especially the (parabolic) induction and restriction functors. 

For any $J \subset \Pi$, let $\CH_J$ be the parabolic subalgebra of $\CH$ generated by $\{\theta_x T_w; x \in X, w \in W_J\}$. This is the affine Hecke algebra for the parabolic subgroup $\tilde W_J=X \rtimes W_J$, for some parameters. 

We may then define the (parabolic) induction and restriction functors:
$$i_J: R(H_{J, \mathbf c'}) \to R(H_{\mathbf c}), \qquad r_J: R(H_{\mathbf c}) \to R(H_{J, \mathbf c'}).$$

\subsection{} For $J \subset \Pi$, we have the central character $\chi_t$ of $H_{J, \mathbf c'}$, here $t$ runs over a complex torus $T^J$ associated to $J$. If $\s \in R(H_{J, \mathbf c'})$, then $\s \chi_t \in R(H_{J, \mathbf c'})$.  Following Bernstein, Deligne and Kazhdan, we say that a form $f \in R(H_{\mathbf c})^*$ is {\it good} if for any $J \subset \Pi$ and $\s \in R(H_{J, \mathbf c'})$, the function $t \mapsto f(i_J(\s \chi_t))$ is a regular function on $T^J$. 

The following result is proved in \cite{BDK} and \cite{Kaz}. 

\begin{thm}
Let $H$ be the Hecke algebra of a $p$-adic group $G$, i.e., the algebra of locally constant, compactly support measures on $G$. Then 

(1) The image of $Tr: \bar H \to R(H)^*$ is the set of good forms. 

(2) The map $Tr: \bar H \to R(H)^*$ is injective. 
\end{thm}

Part (1) is called trace Paley-Wiener Theorem and Part (2) is called the density Theorem. The proofs rely on p-adic analysis. 

%In the rest of this note, we'll generalize these result to extended affine Hecke algebra of generic/arbitrary parameters. It is based on the cocenter-representation duality for extended affine Hecke algebras. 

\section{Cocenter of affine Hecke algebras} 

\subsection{} In this section, we discuss the structure of $\bar \CH$. In order to do this, we first need to generalize Theorem \ref{f-min} to affine Weyl groups. This is rather nontrivial as in general, a conjugacy class in an affine Weyl group contains infinitely many elements, and thus it is impossible to check by computer via brute force whether the desired properties hold even for a very simple group such as $\ZZ^3 \rtimes S_3$. 

\subsection{}
The idea is to use the arithmetic invariants of conjugacy classes in $\tW$.

There are two invariants. The first invariant is given by Kottwitz map $\k: \tW \to \tW/W_a$, sending an element to its coset. The second invariant is given by the Newton point. It is defined as follows. 

For $w \in \tW$, $w^{|W_0|}=t^\l$ for some $\l \in X$. The Newton point $\nu_w$ of $w$ is defined to be $\nu_w=\l/|W_0| \in X_\QQ$. We denote by $\bar \nu_w$ the unique dominant element in $W_0 \cdot \nu_w$. 

The map $$f: \tilde W \to \tilde W/W_a \times X_{\mathbb Q}, \quad w \mapsto (w W_a, \bar \nu_w)$$ is constant on each conjugacy class of $\tilde W$. Each fiber is a union of finitely many conjugacy classes of $\tW$. It allows us to reduce from affine Weyl groups to finite Weyl groups (associated to the Newton point). 

Although not needed in this note, it is worth mentioning that the map $f$ also relates the Frobenius-twisted conjugacy classes of a $p$-adic group with the conjugacy classes of its Iwahori-Weyl group, and plays a crucial role in the study of affine Deligne-Lusztig varieties and Newton strata of Shimura varieties. We refer to \cite{He99} and \cite{He00} for more details. 

\

Now we have the following properties on extended affine Weyl groups, first discovered in \cite{He98} for some classical groups and then proved in \cite{HN14} in general. 

\begin{thm}\label{aff-min}
Let $\tW$ be an extended affine Weyl group and $\co$ be a conjugacy class. Let $\co_{\min}$ be the set of minimal length elements in $\co$. Then 

(1) Each element $w$ of $\co$ can be brought to a minimal length element $w'$ by conjugation by simple reflections which reduce or keep constant the length. 

(2) Any two elements in $\co_{\min}$ are ``strongly conjugate'', i.e., conjugate in the associated Braid group. 
\end{thm}

Similar to the argument for finite Hecke algebras, we have the Iwahori-Matsumoto presentation of the cocenter of extended affine Hecke algebras. 

\begin{thm}
Let $\CH$ be an extended affine Hecke algebra. Then 

(1) For any conjugacy class $\co$ of $\tilde W$ and $w, w' \in \co_{\min}$, the image of $T_w$ and $T_{w'}$ in $\bar \CH$ are the same. We denote it by $T_\co$. 

(2) $\{T_\co\}$ spans of $\bar \CH$. 
\end{thm}

It is much harder to prove that $\{T_\co\}$ is linearly independent. The reason is that both $\bar \CH$ and $R(\CH)$ are of infinite rank and thus one can't compare their rank directly to make the conclusion. For equal parameter case, linearly independence is proved in \cite{HN14} using Lusztig's $J$-ring. The general case is proved in \cite{CH2} using the cocenter-representation duality and will be discussed in \S\ref{6}. 

\subsection{} In order to relate the cocenter to the representations, we also need the Bernstein-Lusztig presentation for $\bar H$. This is based on the following parametrization of conjugacy classes of $\tilde W$ in \cite{HN15}. 

\begin{prop}\label{5.3}
Let $\mathcal A$ be the set of all pairs $(J, C)$, where $J \subset \Pi$ and $C$ is an elliptic conjugacy class of $\tilde W_J$ and $\nu_w$ is dominant for all $w \in C$. Then $$\pi: \mathcal A/\sim \, \xrightarrow{1-1} \text{ conjugacy classes of } \tilde W.$$
\end{prop}

\

Now we have the following presentation of the cocenter. 

\begin{thm}
Let $i_J: \CH_J \to \CH$ be the inclusion and $\bar i_J: \bar \CH_J \to \bar \CH$ the induced map. For any $(J, C) \in \mathcal A$, $$T_\co=\bar i_J(T^J_C).$$
\end{thm}

Note that the relation between Iwahori-Matsumoto basis and Bernstein-Lusztig basis of $\CH$ is complicated and the relation between the minimal length elements in $\co$ (with respect to the length function of $\tW$) and the minimal length elements in $C$ (with respect to the length function of $\tW_J$) is also complicated. It is amazing that these two relations match well. This leads to the above matching between Iwahori-Matsumoto basis and Bernstein-Lusztig basis of the cocenter of $\CH$. 

\section{Elliptic quotient and Rigid quotient}\label{6}

\subsection{} Elliptic representation theory, introduced by Arthur \cite{A}, studies
the Grothendieck group of certain representations of a Lie-theoretic
group modulo those induced from proper parabolic subgroups. The elliptic
theory of representations of semisimple $p$-adic groups and
Iwahori-Hecke algebras was further studied intensively, e.g., Schneider-Stuhler
\cite{SS}, Bezrukavnikov \cite{Bez}, Reeder \cite{Re}, Opdam-Solleveld
\cite{OS}. 

For an affine Hecke algebra $H_{\mathbf c}$, the {\it elliptic quotient} is defined to be $$R(H_{\mathbf c})_{ell}=R(H_{\mathbf c})_\CC/\sum_{J \subsetneqq \Pi} i_J(R(H_{J, \mathbf c'})_\CC).$$

Opdam and Solleveld in \cite{OS} studied the affine Hecke algebras for positive parameters and showed that 

\begin{thm}
Let $\mathbf c$ be a positive parameter function on $H$. Then  

(1) The dimension of $R(H_{\mathbf c})_{ell}$ is at most the number of elliptic conjugacy classes of $\tW$.

(2) The inequivalent discrete series forms an orthogonal set of $R(H_{\mathbf c})_{ell}$. 
\end{thm}

In particular, one has an upper bound for the number of irreducible discrete series for affine Hecke algebra of positive parameters. A lower bound can be obtained by counting the central characters of the discrete series. This leads to the classification of irreducible discrete series for affine Hecke algebras of positive parameters in \cite{OS2}. 

\subsection{} The method of Opdam-Solleveld is analytic, passing from $H$ to its Schwartz algebra, a certain topological completion of $H$. For Hecke algebras of $p$-adic groups, the elliptic quotient and its relation with the trace map was also studied in \cite{BDK} and analytic methods were used to obtain an upper bound of the dimension of the elliptic quotient. 

\subsection{} Here is our motivation to develop a different method to study elliptic representation theory for affine Hecke algebras. 

First, we would like to understand the affine Hecke algebras for arbitrary parameters (not just the positive parameters) and for representations over a field of positive characteristics. Thus we'd prefer to have a more algebraic method. 

Second, we'd like to put the elliptic quotient in the framework of ``cocenter-representation duality''. Namely, the elliptic quotient $R(H_{\mathbf c})_{ell}$ corresponds to a subspace of $\bar H$ via the trace map $Tr: \bar H \times R(H)_\CC \to \CC$. What is this subspace? Results in \cite{BDK} indicates that this subspace is very complicated and may not have a nice explicit description. 

\subsection{}\label{rig} The idea in \cite{CH2} is to replace the elliptic quotient by the so-called rigid quotient. For simplicity, we only give the definition for extended affine Hecke algebras associated to semisimple root data. We simply write $H$ for $H_{\mathbf c}$.

We define a subspace of $\bar H$ by the arithmetic invariants: $$\bar H_{\Pi} =\{T_\co; \nu_\co=0\}.$$ and a quotient space of $R(H)_\CC$ by $$R(H)_{rig}=R(H)_\CC/\<i_J(\s)-i_J(\s \chi_t); J \subset \Pi, \s \in R(H_{J}), t \in T^J\>,$$ the quotient of $R(H)$ be the difference of induced modules. We call $R(H)_{rig}$ the rigid quotient of $H$. Both $\bar H_1$ and $R(H)_{rig}$ are finite dimensional. 

The main result in \cite{CH2} is that

\begin{thm}\label{rig-pairing}
(1) For generic parameters, the trace map $Tr: \bar H \times R(H)_\CC \to \CC$ induces a perfect pairing \[\tag{a} Tr: \bar H_{\Pi}  \times R(H)_{rig} \to \CC.\]

In particular, $\dim R(H)_{rig}=\dim \bar H_{\Pi} $ equals the number of conjugacy classes in $\tW$ whose Newton points equal to zero. 

(2) For arbitrary parameter, the induced map $Tr: \bar H_{\Pi}  \to R(H)_{rig}^*$ is surjective. 
\end{thm}

We call Theorem \ref{rig-pairing} (a) the ``cocenter-representation'' duality for affine Hecke algegbras. Now we discuss some consequences. 

\begin{thm}
For generic parameters, 

(1) $\{T_\co\}$ is a basis of $\bar H$. 

(2) The image of the map $Tr: \bar H \to R(H)^*$ is $R(H)^*_{good}$. 

(3) The map $Tr: \bar H \to R(H)^*$ is injective. 
\end{thm}

\subsection{} In fact, we have a natural projection map $R(H)_{rig} \to R(H)_{ell}$. We may choose a basis of $R(H)_{rig}$, which includes a basis of $R(H)_{ell}$ as a subset. We also have a standard basis on $\bar H_{\Pi} $, the Iwahori-Matsumoto/Bernstein-Lusztig basis. The matrix associated to $Tr: \bar H_{\Pi}  \times R(H)_{rig} \to \CC$ is in general, far from being block triangular. This is the reason that the ``dual'' of elliptic quotient is hard to understand and the reason that we consider the rigid quotient instead. 

On the other hand, one gets the upper bound of the dimension of the elliptic quotient space from the estimate on the rigid quotient space. 

\begin{cor}
For arbitrary parameters, the dimension of $R(H_{\mathbf c})_{ell}$ is at most the number of elliptic conjugacy classes of $\tW$.
\end{cor}

We also have the following deformation theorem \cite{CH2}. 

\begin{thm}
For a generic parameter function $\mathbf c$, there exists a basis $\{V_{\pi, \mathbf c}\}$ of $R(H)_\QQ$ such that for any $w \in \tW$ and any $\pi$, the action of $T_w$ on $V_{\pi, \mathbf c}$ depends analytically on the parameter function $\mathbf c$. 
\end{thm}

For affine Hecke algebra with positive paramenters, a similar result was obtained by Opdam and Solleveld in \cite{OS2} using Schwartz algebras. The idea of our proof (without the restriction on positive parameters) is to use Lusztig's graded affine Hecke algebras \cite{L1} and the deformation theorem for graded affine Hecke algebras obtained in \cite{CH}. The passage from affine Hecke algebras to graded affine Hecke algebras is analytic. This is the reason that we have the analytic deformation theorem here. However, we expect that there exists a family of representations depending algebraically on the parameter function. 

\section{The "Modular case"}

\subsection{} Now we move to the more ``modular case''. Here instead of $\CC$, we consider an algebraically closed field $\kk$ of positive characteristic. We also consider the case where the parameter is a root of unity. We have the extended affine Hecke algebra $H_{\mathbf c, \kk}$. The situation is much more complicated and one can't expect to have a perfect pairing as in Theorem \ref{rig-pairing} (1). However, we expect that the cocenter-representation duality still holds (under some modification) and provide useful information on the representations of affine Hecke algebras in the ``modular case''. 

\subsection{} %First by a deformation argument, the basis theorem and the trace Paley-Wiener Theorem holds for arbitrary $\mathbf c$ and arbitrary $\kk$ and the filtration can be defined in the same way. 

We still have the map $Tr: (\bar H_{\mathbf c, \kk})_\Pi \to R(H_{\mathbf c, \kk})^*_{rig}$, but it fails to be injective. Since the parahoric subalgebra $H_J^{fin}$ may not be semisimple, thus $T_\co$ for some $\co \in W_J$ may lie in the kernel of the map $Tr: (\bar H_{\mathbf c, \kk})_\Pi \to R(H_{\mathbf c, \kk})^*_{rig}$. It is interesting to see whether the whole kernel comes from the non-semisimplicity of the parahoric subalgebras. This leads to the following (conjectural) estimate on the dimension of the elliptic quotient. 

\subsection{} Recall that $\tW=W_a \rtimes \Omega$. For any $J \subset \tS$ with $W_J$ finite, let $\Omega_J$ be the subgroup of $\Omega$ that stabilizes $J$ and $W_J^\sharp=W_J \rtimes \Omega_J$. Let $\ci$ be the set of all $J$ with $W_J^\sharp$ maximal. Then we ask 

\begin{ques}\label{8.1}
Is $\text{rank}(\bar R(H_{\mathbf c, \kk})_{ell})=\sum_{J \in \ci/\sim} \text{rank}(\bar R(H^{fin}_{J, \mathbf c, \kk})_{ell})$? 
\end{ques}

\subsection{} Now we provide some examples. 

First we consider the case where $\tW$ is the Iwahori-Weyl group of $SL_3$. Then $\Omega=\{1\}$ and $\ci$ consists of three subsets of $\tS$: $\{0, 1\}, \{1, 2\}, \{0, 2\}$. The rank is given as follows. 

%\begin{table}
\[\begin{tabular}{|c | c | c|}
\hline
$SL_3$ &  $char(\kk) \neq 3$ & $char(\kk)=3$ \\ \hline
$\Phi_3(q) \neq 0$ &    3       & 3 \\ \hline
$\Phi_3(q)=0$       &    0       & 0 \\
\hline
\end{tabular}\]
%\end{table}

Then we consider the case where $\tW$ is the Iwahori-Weyl group of $PGL_3$. Then $\Omega=\ZZ/3 \ZZ$ and $\ci$ consists of four subsets of $\tS$: $\{0, 1\}, \{1, 2\}, \{0, 2\}$ and $\emptyset$. The first three subset are conjugate in $\tW$ and $W_\emptyset^\sharp=\Omega$. The rank is given as follows. 

\[
\begin{tabular}{|c | c | c|}
\hline
$PGL_3$ &  $char(\kk) \neq 3$ & $char(\kk)=3$ \\ \hline
$\Phi_2(q) \neq 0$ &    3       & 1 \\ \hline
$\Phi_2(q)=0$       &     2       & 0 \\
\hline
\end{tabular}
\]

\subsection{} Let us discuss the case where $\kk=\CC$ and $\mathbf c$ is a constant (i.e. equal parameter case). In this case, $H^{fin}_{J, \mathbf c, \CC}$ is semisimple for any $J$ if and only if the parameter is not a root of Poincar\'e polynomial for $W_0$. The conjectural equality in Question \ref{8.1} for $H_{\mathbf c, \CC}$ and its parabolic subalgebras indicates how the ``number'' of irreducible representations changes when the parameter changes to a root of Poincar\'e polynomial. 

In particular, if the parameter is not a root of Poincar\'e polynomial, then the parametrization of irreducible representations is independent of the choice of the parameter. For a root datum $\Phi$ whose associated group has simply connected derived group, this is proved by Xi in \cite{Xi}.

\section*{Acknowledgment} We thank D. Ciubotaru, G. Lusztig, J. Michel, S. Nie and E. Opdam for useful discussions.

%\section{Relation to $p$-adic groups (?)} 

%\section{Other cases} 

%graded affine Hecke algebras, $0$-Hecke algebras, DAHA(?)

\end{document}